\documentclass[aps,preprint,superscriptaddress,floatfix,amsmath,amssymb]{revtex4}

\usepackage{inputenc}
\usepackage[a-1b]{pdfx}
\usepackage{amsmath}
\usepackage{multirow}
\usepackage{booktabs}
\usepackage{amssymb}
\usepackage{amsthm}
\usepackage{algcompatible}
\usepackage{algpseudocode}
\usepackage{newfloat}
\usepackage{bm}
\usepackage{graphicx}
\usepackage{subfig}
\usepackage{hyperref}
\usepackage{float}
\usepackage{subfig}
\usepackage{array}

\begin{document}

\title{An Initial Condition-Dependent Neural Network Approach for Optimal Control Problems}

\author{Mominul Rubel}
\email{mrubel@mst.edu}
\affiliation{Engineering Management and Systems Engineering, Missouri University of Science and Technology, Rolla, MO}

\author{Gabriel Nicolosi}
\email{gabrielnicolosi@mst.edu}
\affiliation{Engineering Management and Systems Engineering, Missouri University of Science and Technology, Rolla, MO}

\date{Februrary 10, 2025}

\begin{abstract}
{\small In this work, we investigate an indirect approach for the numerical solution of optimal control
problems via neural networks. A customized neural network is constructed,
where optimal state, co-state and control trajectories are approximated by minimizing the
underlying parameterized Hamiltonian, relying on Pontryagin's Minimum Principle. Departing
from previous results reported in the literature, we propose novel, modified networks with
both time and trajectory initial condition as inputs. Numerical results demonstrate the ability of neural networks to integrate both time and initial condition information in solving optimal control problems. Finally, it is empirically demonstrated that approximation accuracy may be enhanced through a structural modification incorporating an intermediate layer of Fourier coefficients.}
\end{abstract}

\maketitle

\section*{Keywords}
Optimal control problem, Pontryagin's Minimum Principle, Neural networks, Fourier series

\section{Introduction}
Optimal control problems (OCPs) are optimization problems that seek to find a control input that guides the dynamics of a system while seeking to minimize  or maximize a given performance index represented by an objective function. In the context of management sciences and economics, a standard treatment of the field can be found in \cite{sethi2022optimal}. Traditional methods, such as Pontryagin’s Minimum Principle (PMP) (necessary conditions) and Hamilton-Jacobi-Bellman (HJB) equations (sufficient conditions), have been commonly applied to solve such problems. However, for systems of real, practical importance, typically describing complex, nonlinear behavior, these classical approaches are rather analytically limited and intractable, often failing to provide a closed-form solution. Thus, the solution of many OCPs must rely on the efficient application of numerical techniques. 

Following the mass dissemination and use of machine learning models in the sciences and engineering, we consider the use of artificial neural networks (NNs) as surrogates for optimal control laws (and their induced optimal trajectories), offering a viable alternative to the current apparatus of numerical optimal control. Departing from the work developed in \cite{EP13, NFG22, NG23, NFG25}, we generalize and integrate their approach, considering a family of optimal control problems over a given range of initial conditions. The main contribution of this present work lies in extending and integrating these models to accommodate a set of initial conditions, thus approximating optimal control and trajectory surfaces with a single set of parameters: weights, biases and possibly Fourier coefficients.

This paper is structured as follows: Section~(\ref{sec:literature}) presents a brief overview of the main frameworks in numerical optimal control and recent contributions to the field from a machine learning perspective. Section~(\ref{Sec: Problem Construction}) constructs two different numerical methods building upon and adding to others found in the recent literature. More specifically, two neural network-based approaches are presented to approximate the solution of OCPs over a range of system initial conditions: a traditional shallow neural network and a neural network equipped with a layer of Fourier coefficients, with the latter exploiting the approximation capabilities of both neural networks and Fourier series. Section~(\ref{sec:computation}) presents numerical results for three different OCPs used as testbeds to methods presented herein. Finally, Section~(\ref{sec:conclusion}) discusses the implications of these findings and draws a future line of work based on the results obtained so far.

\section{Related Literature}\label{sec:literature}
The field of numerical optimal control is vast, with applications ranging from distinct fields such aeronautical engineering to economics. For an overview, the reader is pointed to \cite{Practical2020}. The methods are typically classified into two: direct and indirect. In the former, the target functions (control, state and co-state) are parameterized by approximating functions, such as orthogonal basis functions (\textit{e.g.}, Chebyshev, Legendre or Fourier). In this way, the original OCP is transformed into a finite dimensional optimization problem, where the dynamics to be respected is priced out into the objective function and standard nonlinear programming algorithms utilized to find a set of optimal parameters that best approximate the target functions. This approach is typically subdivided into shooting and collocation methods \cite{K17}. In its turn, indirect methods take into consideration necessary conditions (sufficient under convexity assumptions), namely the PMP. Even though this approach also involves target function parametrization, it indirectly minimizes the performance index by minimizing the so-called Hamiltonian function with respect to the control input, thus justifying its name. This is briefly introduced in Section~\ref{Sec: Problem Construction} and it is the approach considered in this work.   

Typically, a numerical solution of a given OCP is found for a given fixed initial condition. Such is the case found in \cite{EP13} where, instead of recurring to traditional approximation schemes, suggests an NN approximation for solving OCPs within an indirect framework, thus incorporating a parameterize Hamiltonian into the NN loss function. Similarly and more recently, \cite{BAA22} proposes an indirect method for solving OCPs based on the framework of neural differential equations \cite{CRBD18}. \cite{NFG22, NG23, NFG25} first proposed the (direct) approximation of an OCP for arbitrary initial conditions within a given range. In these series of papers, multidimensional Fourier series are used as approximators for a family of OCPs, generating time and initial-condition dependent optimal control and state surfaces.             

Next, we combine the frameworks presented above into two similar but distinct methods: (i) a time and initial condition-dependent NN, generalizing the method suggested in \cite{EP13} for a range of initial conditions (thus, doing so inspired by the ideas in \cite{NFG22, NG23, NFG25}); (ii) an initial condition-based NN with a Fourier layer. Both approaches seek to approximate a family of OCPs starting from different initial conditions with a single set of parameters.

\section{Problem Construction}\label{Sec: Problem Construction}
Consider the following OCP of Bolza form, \textit{i.e.}, the combination of the Lagrange and Mayer forms
\begin{equation}
\begin{split}
    &\min_{u} \int_{0}^{T} f(x(t), u(t), t) \, dt + \Psi(x(T), T)\\
    &\text{\textit{s.t.}} \quad \dot{x}(t) = g(x(t), u(t), t) \\
    & \quad\quad x(0) = x_0,
\end{split}
\label{eq:OCP1}
\end{equation}

where the functions \( x: [0, T] \to \mathbb{R} \) and \( u : [0, T] \to \mathbb{R} \) describe the state of the system and the control input, respectively, in the time interval \( [0, T]\). The running cost function \(f(x, u, t)\) and the system dynamics \(g(x, u, t)\) are assumed to be continuously differentiable. Suppose \(X_0 = \{x_0^1, x_0^2, \dots, x_0^l\}\) is a finite ordered set of initial states. The goal is to determine the control trajectory \(u(t,x_0)\) that solves OCP~(\ref{eq:OCP1}) for all $x_0$ belonging to the convex hull of $X_0$, which we shall denote by $\mathcal{X}_0$. Consider the Hamiltonian function for a fixed \(x_0\), defined by
\begin{equation}
    H(x, u, \lambda, t) :=  f(x, u, t) + \lambda  g(x, u, t),
    \label{eq:hamiltonian}
\end{equation}

where \( \lambda : [0, T] \to \mathbb{R} \) is the so-called costate function. PMP states that a solution to OCP~(\ref{eq:OCP1}), \textit{i.e.}, an optimal control $u^*(t)$, its induced trajectory $x^*(t)$ and costate $\lambda^*(t)$ are, necessarily, solutions to the following system of differential equations
\begin{equation}
    \frac{\partial H}{\partial x} = -\dot{\lambda}(t), \quad 
    \frac{\partial H}{\partial \lambda} = \dot{x}(t), \quad 
    \frac{\partial H}{\partial u} = 0.
    \label{eq:pmps}
\end{equation}

Next, we introduce two NN-based approaches that incorporate the conditions in Equation~(\ref{eq:pmps}) into the loss function, as initially done in \cite{EP13}, but now taking into account the set $\mathcal{X}_0$, thus aiming to obtain numerical approximations for a family of OCPs defined for all initial conditions in $\mathcal{X}_0$.

\subsection{Time and Initial Condition-Dependent NN}\label{subsec: NN}
In this method, a customized NN is composed of three parallel, sub-networks, each designed to approximate the (optimal) target functions state \(x^*(t)\), costate \(\lambda^*(t)\), and control \(u^*(t)\). This NN takes both time \(t \in [0,T]\) and initial condition \(x_0 \in X_0\) as inputs and generates three respective outputs, each corresponding to one target function. By passing these output signals to the PMP conditions embedded in the loss function, these approximations are forced to simultaneously solve the underlying OCP for all initial conditions in $X_0$. Boundary conditions are satisfied by designing so-called trial solutions, shown in Equation~(\ref{eq:trialsolns}). The NN outputs for the state approximation is defined as

\begin{equation}
        n_x(t, x_0) = \boldsymbol{v_x}^\top \sigma(\boldsymbol{z_x}),\quad  \boldsymbol{z_x} = \boldsymbol{w_x}^\top [t\quad x_0]^\top + \boldsymbol{b_x},
    \label{eq: nns}
\end{equation}

where \(\sigma(\cdot)\) is the sigmoid activation function, \(\boldsymbol{v_x}\) is a vector of weights of the output layer, \(\boldsymbol{w_x}\) a matrix of weights, and \(\boldsymbol{b_x}\) a vector containing the biases of the NN hidden layer. $n_u$ and $n_\lambda$ are defined likewise, with their own corresponding parameters. We pack all of NN parameters into a single vector, \(\boldsymbol{\Phi}\). The trial solutions are defined as

\begin{equation}
        \widehat{x}(t, x_0) = x_0 + t n_x, \quad
        \widehat{\lambda}(t, x_0) = (t-T) n_\lambda\quad \text{and} \quad
        \widehat{u}(t, x_0) = n_u, 
        \label{eq:trialsolns}
\end{equation}

and ensure that \(\widehat{x}(0) = x_0\), \(\widehat{\lambda}(T) = 0\) (transversality condition). To ensure that the trial solutions satisfy PMP conditions, the trial solutions are substituted in Equation~(\ref{eq:hamiltonian}), yielding the trial Hamiltonian, emphasizing the role of the initial condition by adding a $x_0$ as an argument to the original formulation introduced in \cite{EP13}. The trial Hamiltonian is written as

\begin{equation}
\widehat{H}(\widehat{x}, \widehat{\lambda}, \widehat{u}, t, x_0) = f(\widehat{x}, \widehat{u}, t) +\widehat{\lambda}g(\widehat{x}, \widehat{u}, t).    
\end{equation}

The derivatives of the trial Hamiltonian are then used to define error functions corresponding to the PMP conditions

\begin{equation}
        E_1(t, x_0; \boldsymbol{\Phi}) = \frac{\partial \widehat{H}}{\partial \widehat{x}} + \dot{\widehat{\lambda}}, \quad
        E_2(t, x_0; \boldsymbol{\Phi}) = \frac{\partial \widehat{H}}{\partial \widehat{\lambda}} - \dot{\widehat{x}}, \quad
        E_3(t, x_0; \boldsymbol{\Phi}) = \frac{\partial \widehat{H}}{\partial \widehat{u}}
        \label{eq:trialham}
\end{equation}
composing the loss function
\begin{equation}
    E(t, x_0; \boldsymbol{\Phi}) = E_1^2 + E_2^2 + E_3^2.
    \label{eq:lossfun}
\end{equation}
 Equation~(\ref{eq:lossfun}) is used to assess the solution accuracy by measuring the deviation from the necessary PMP conditions and we seek to minimize it over a grid of points defined by the Cartesian product \(\mathcal{T} \times X_0\), where \(\mathcal{T}\) is a set of discrete time points over the interval \([0, T]\). The resulting unconstrained optimization problem is given by

\begin{equation}
    \min_{\boldsymbol{\Phi}} \textstyle\sum_{\;x_0\in X_0} \textstyle\sum_{\;t \in \mathcal{T}} E(t, x_0; \boldsymbol{\Phi}).
    \label{eq:unconstrainedopt}
\end{equation}

We denote by \(\boldsymbol{\Phi^*}\) a solution to the optimization problem in Equation~(\ref{eq:unconstrainedopt}), and it represents the best set of parameters that satisfies Equation~(\ref{eq:pmps}), thus parameterizing all of the target functions for all initial conditions in $\mathcal{X}_0$. 

\subsection{Initial Condition-Dependent NN with a Fourier Layer}\label{subsec: NNF}
\begin{figure}[ht]
    \centering
    \includegraphics[width=0.6\linewidth]{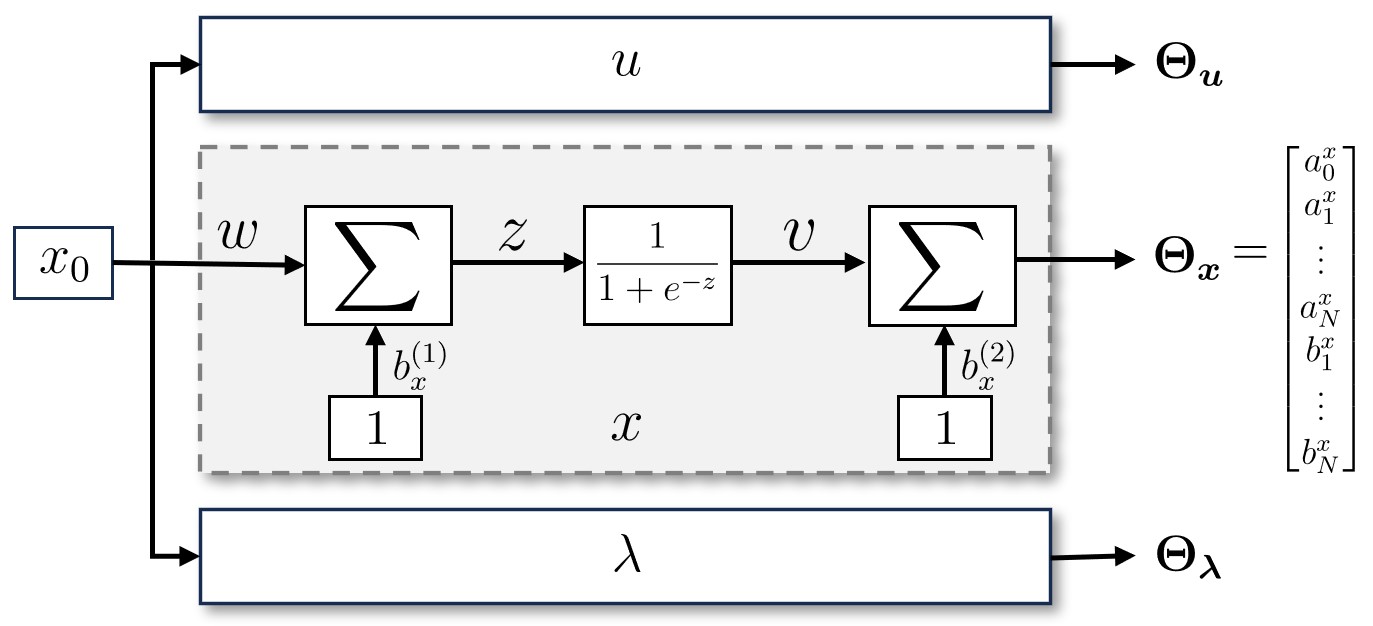}
    \caption{Schematic representation of the initial condition-dependent NN with a Fourier layer. For illustrative purposes, the network structure is shown only for the state approximation but it is identical for both control and co-state functions.}
    \label{fig:theNN}
\end{figure}

We now consider only the initial condition \(x_0\) as the input of each of the three sub-networks, as shown in Equation~(\ref{fig:theNN}). Additionally, an output layer of Fourier coefficients for each target function is added to the NN. The idea behind such design spans from the work in \cite{NFG22, NG23, NFG25}. As such, the trial functions giving rise to the desired approximations are now parameterized by distinct finite Fourier series, obtained through a map between an initial condition $x_0$ and a finite set of Fourier coefficients, thus aiming to leverage the approximation capabilities of both NNs and Fourier decompositions for the numerical solution of OCPs. 

The vector of Fourier coefficients $\boldsymbol{\Theta}_x(x_0)$, parameterizing the approximate optimal state trajectory for a given $x_0 \in X_0$, similarly to Equation~(\ref{eq: nns}), is defined as
\begin{equation}
        \boldsymbol{\Theta}_x(x_0) = \boldsymbol{v_x}^\top \sigma(\boldsymbol{z_x}) + \boldsymbol{b^{(2)}_x},\quad  \boldsymbol{z_x} = \boldsymbol{w_x} x_0 + \boldsymbol{b^{(1)}_x}, 
\end{equation}

with $\boldsymbol{\Theta}_u(x_0)$ and $\boldsymbol{\Theta}_\lambda(x_0)$ defined accordingly. Furthermore, the trial functions of Equation~(\ref{eq:trialsolns}) become finite Fourier series, now being written as
\begin{equation}
    \begin{aligned}
        \widehat{u}(t, x_0; \boldsymbol{\Theta_u}) &:= a^u_0 + \textstyle\sum_{m=1}^M \left( a^u_m \cdot \sin\left(\tfrac{m \pi t}{T}\right) + b^u_m \cdot \cos\left(\tfrac{m \pi t}{T}\right) \right), \\
        \widehat{x}(t, x_0; \boldsymbol{\Theta_x}) &:= \left(x_0 - \textstyle\sum_{n=1}^N b^x_n\right) + \textstyle\sum_{n=1}^N \left( a^x_n \cdot \sin\left(\tfrac{n \pi t}{T}\right) + b^x_n \cdot \cos\left(\tfrac{n \pi t}{T}\right) \right), \\
        \widehat{\lambda}(t, x_0; \boldsymbol{\Theta_\lambda}) &:= \left(\lambda_T - \textstyle\sum_{n=1}^N b^\lambda_n \cdot (-1)^n\right)+ \textstyle\sum_{n=1}^N \left( a^\lambda_n \cdot \sin\left(\tfrac{n \pi t}{T}\right) + b^\lambda_n \cdot \cos\left(\tfrac{n \pi t}{T}\right) \right).
    \end{aligned}
\end{equation}

We note that we have written these Fourier series with slight modifications to their constant terms so that they satisfy boundary conditions such as \(\widehat{x}(0, x_0; \boldsymbol{\Theta_x}) = x_0\) and \(\widehat{\lambda}(T, x_0; \boldsymbol{\Theta_\lambda}) = \lambda_T\) when appropriate. Substituting \(\widehat{x}, \widehat{\lambda},\) and \( \widehat{u}\) in Equation~(\ref{eq:hamiltonian}), the loss function and the unconstrained optimization problem are constructed in the same manner as described in Section~(\ref{subsec: NN}). The solution to this problem provides a map between initial condition and a set of Fourier coefficients used for the approximation of the desired target functions for all $x_0 \in \mathcal{X}_0$. 

\section{Computational Results} \label{sec:computation}
\begin{figure}[htpb]
   \centering
   \includegraphics[width=0.60\textwidth]{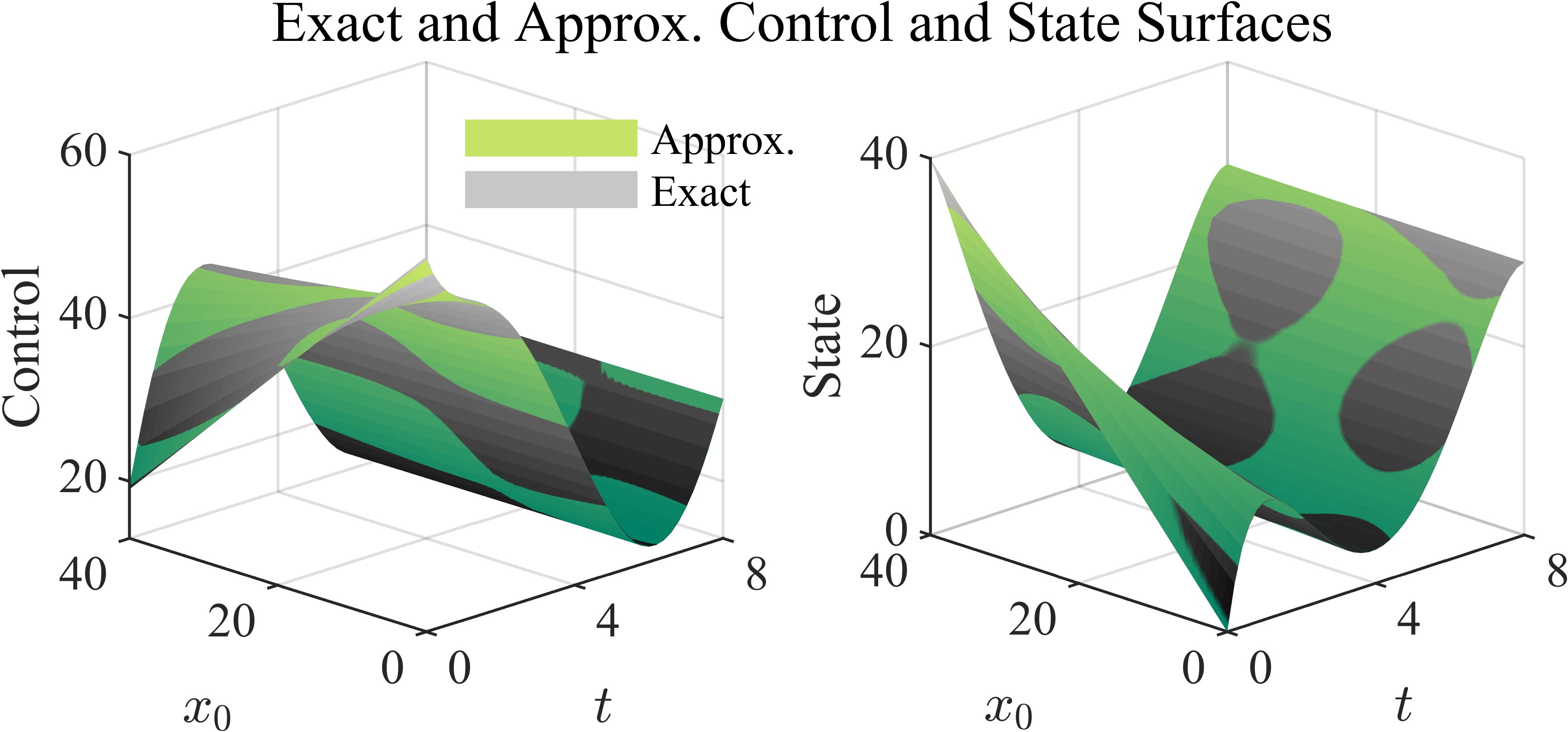}
   \qquad
   \includegraphics[width=0.3\textwidth]{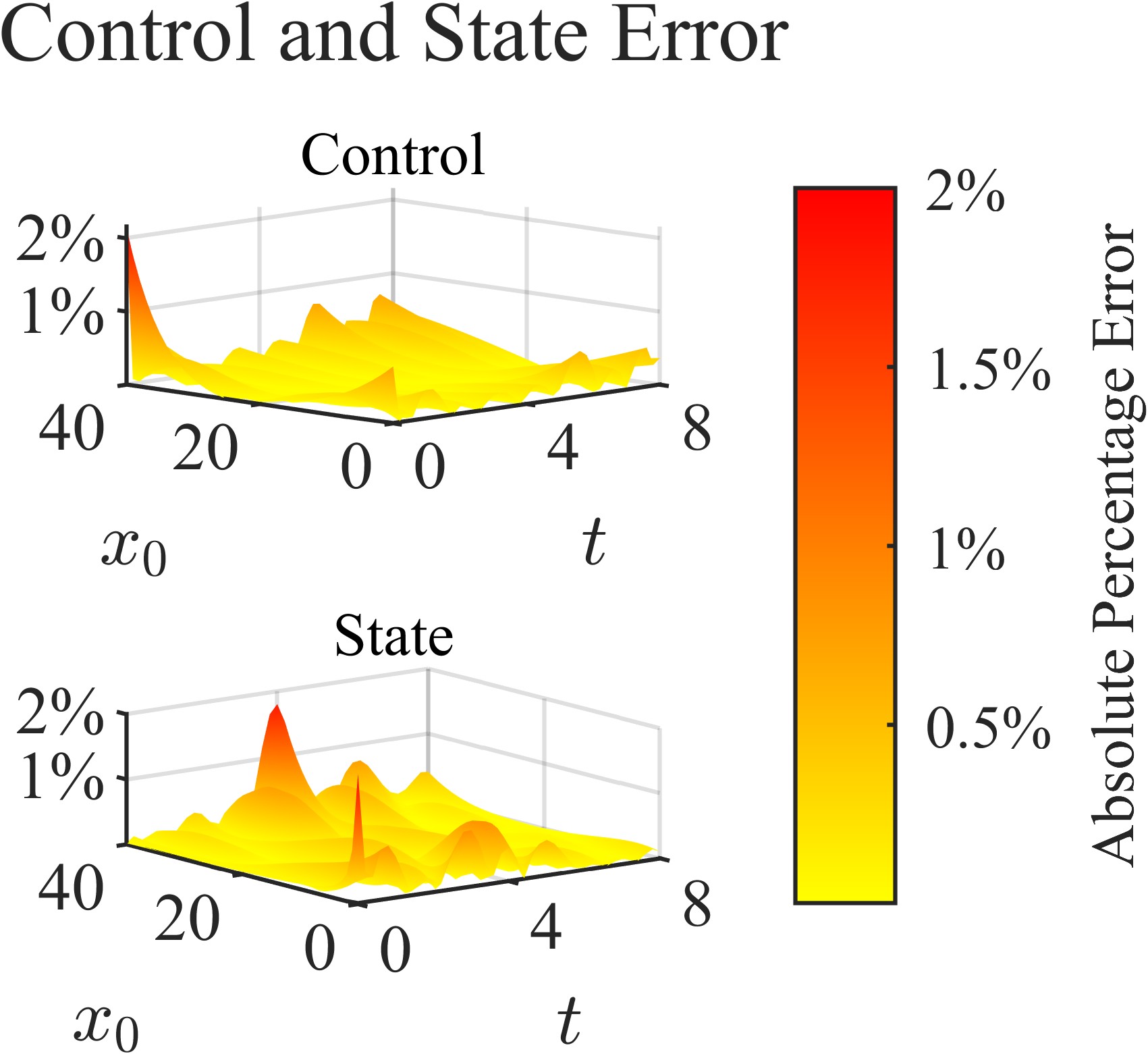}
   \caption{Results for OCP 3 for $M=N=5$ and $I = 6$ (with the Fourier layer). (Left) Approximated and the exact control and state surfaces, respectively. (Right) Absolute Percentage Error \((APE)\) in the control and state approximations, respectively, \(\forall(t, x_0) \in \mathcal{T} \times X_0\).} \label{fig:mainfig}
\end{figure}

\begin{figure}[ht!]
    \centering
    \includegraphics[width= 0.9\linewidth]{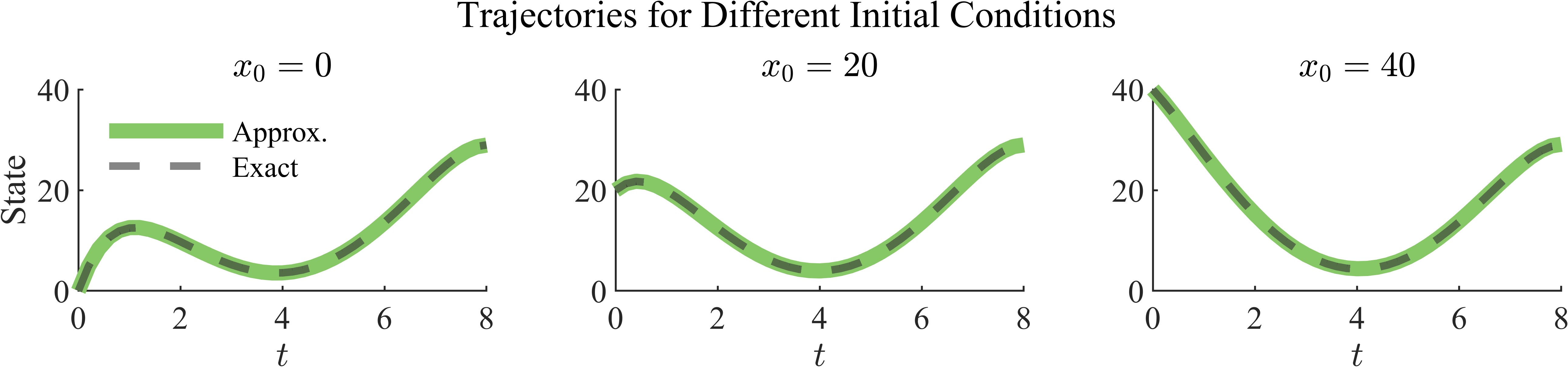}
    \caption{Sample trajectories of the state surfaces shown in Figure~(\ref{fig:mainfig}) for three distinct initial conditions.}
    \label{fig:trajectories}
\end{figure}

We provide numerical results for the three OCPs defined in Table~(\ref{tab: OCPS}) for the methods presented. Plots are shown in Figs.~(\ref{fig:mainfig}) and ~(\ref{fig:trajectories}) for OCP 3 for the method considered in Section~(\ref{subsec: NNF}), which describes an inventory control problem. This problem is taken from \cite{sethi2022optimal} and the goal is to keep inventory and production optimally close to their target (safety stock) level given an initial inventory level $x_0$. The parameters used were \(\rho = 0\), \(h=c=1\), \(x_{target} = 15\), and \(u_{target}=30\). OCP 1 and 2 in Table~(\ref{tab:summary_table}) are extracted from \cite{EP13}. 
 
\begin{table}[ht]
\centering
\renewcommand{\arraystretch}{1.2}
\setlength{\abovedisplayskip}{0pt}
\setlength{\belowdisplayskip}{0pt}
\caption{OCPs and their respective time and initial condition domains for the experiments reported in Table~(\ref{tab:summary_table}).}
\scriptsize{
\begin{tabular}{| >{\centering\arraybackslash}m{3.7cm} | >{\centering\arraybackslash}m{4cm} | >{\centering\arraybackslash}m{7.5cm} |}
\hline
\textbf{OCP 1 (Lagrange)} & \textbf{OCP 2 (Mayer)} & \textbf{OCP 3 (Inventory Control -- Lagrange)}\\ \hline\hline
    \begin{equation*}
        \begin{aligned}
            &\min_{u}\int_0^T \left[x^2(t) + u^2(t)\right] \, dt \\
            &\text{\textit{s.t.}}\quad \dot{x}(t) = u(t)\\
            &\quad\quad x(0) = x_0
        \end{aligned}
    \end{equation*}& 
    \begin{equation*}
    \begin{aligned}
        &\min_{u} \quad -x(2) \\
        &\text{\textit{s.t.}} \quad \dot{x}(t) = \frac{5}{2}(-x + ux - u^2)\\
        &\quad\quad x(0) = x_0
    \end{aligned}
    \end{equation*} & 
    \begin{equation*}
    \begin{aligned}
        &\min_{u} \int_{0}^{T} e^{\rho t} \left[\frac{h}{2}(x-x_{target})^2 + \frac{c}{2}(u-u_{target})^2\right] \, dt\\
        &\text{\textit{s.t.}} \quad \dot{x}(t) = u(t) - S(t); \quad S(t) = t^3- 12t^2+ 32t+ 30\\
        &\quad\quad x(0) = x_0
    \end{aligned}
    \end{equation*}\\ \hline
    \(t \in [0,1]\), \( X_0 = \{0.00, 0.05, \dots, 1\} \)&
    \(t \in [0,2]\), \( X_0 = \{0.00, 0.05, \dots, 1\} \) &
    \(t \in [0,8]\), \( X_0 = \{0, 1, \dots, 40\} \) \\ \hline
\end{tabular}
}
\label{tab: OCPS}
\end{table}

Computational results involving the application of both methods suggested in this paper, for the solution of the OCPs in Table~(\ref{tab: OCPS}) are summarized in Table~(\ref{tab:summary_table}). We shall now denote the methods of Sections~(\ref{subsec: NN}) and (\ref{subsec: NNF}) as Method 1 and 2, respectively. To evaluate the performance of our proposed methods, in Table~(\ref{tab:summary_table}), we report: root mean square error ($RMSE_u$), mean absolute error ($MAE_u$), and mean absolute percentage error ($MAPE_u$) with respect to the control approximation. $J^*_{\%error}$ denotes the average percentage of error for the optimal objective function value. All these error measures are averaged over the grid defined by $\mathcal{T} \times X_0$. The reader should note that the first column of Table~(\ref{tab:summary_table}) consists of the two indirect methods discussed in this paper, and, for the sake of comparison, the (direct) method as seen in \cite{NFG25} is also reported for OCP 3. The columns labeled $OCP$, $M$, $N$ and $I$ contain, respectively: the OCP, number of Fourier terms in the control approximation, number of Fourier terms in state approximation and the number of neurons in the hidden layer. All error measures are reported for a training set ($X_0$, defined in Table~(\ref{tab: OCPS})) and a testing set, containing a larger discrete set of points belonging to $\mathcal{X}_0$ that do not belong to $X_0$.

Both methods are capable of approximating the target functions for all three OCPs, attaining \(RMSE_u\) and \(MAE_u\) in the order of $10^{-4}$ in some cases. \(MAPE\) can become relatively large when the actual values being approximated are closer to zero, as observed, for example, in Exp.\ 1 and most of the experiments for OCP 2. As expected, for Method 1, the more neurons in the hidden layer, the better the approximations become (compare Exp.\ 3 and 6).  Furthermore, Exp.\ 6 ($I=30$) and Exp.\ 21 ($I=6$) demonstrate that the addition of a Fourier layer for solving OCP 3 is capable of attaining the same error performance, however with many less neurons in the hidden layer, thus, becoming the most relevant experimental result of this paper. All the experiments related to OCP 1 attained similar performance (Exps.\ 7--12), suggesting no clear relation between $M$, $N$, and $I$. The results obtained from the direct method (Exps.\ 26--28) are given at the end of Table~\ref{tab:summary_table}. With smaller values of $M$ and $N$, the NN-based approaches perform better when compared to the direct method. This result is expected since direct methods completely disregard the existence of the PMP conditions.

\section{Conclusion and Future Works}\label{sec:conclusion}
By incorporating the ideas of \cite{NFG22, NG23, NFG25} into the framework initially proposed by \cite{EP13}, this paper confirmed that neural networks can effectively be integrated into the numerical apparatus to solve OCPs. A single hidden layer and a few neurons were sufficient to provide accurate approximations. Furthermore, as reported experimentally, the inclusion of a Fourier layer may reduce the required number of neurons in the sub-networks hidden layers. The apparent advantage of including a Fourier is that the interpolation for the approximated values is implicit by the Fourier series, instead of linear interpolation implicit in Method 1. This structural addition was shown to be viable without compromising the performance of the underlying approximations.  Future work will investigate the performance of these methods for systems with larger state spaces and attempt to formalize aspects of the computational complexity of both methods, thus far only observed empirically.   

\begin{table}[H]
    \caption{Summary of computational results.}\label{tab:summary_table}
    \centering
    \renewcommand{\arraystretch}{0.7}
    \setlength{\tabcolsep}{1pt}
    \footnotesize{
    \begin{tabular}{|c|c|c|c|c|c|c|c|c|c||c|c|c|c|}
        \hline
        \textbf{} & \textbf{} & \textbf{} & \textbf{} & \textbf{} & \textbf{} & \multicolumn{4}{c||}{\textbf{Training}} & \multicolumn{4}{c|}{\textbf{Testing}} \\
        \cline{7-14}
        & \textit{\textbf{Exp.}} & \textit{OCP} & \textbf{\(M\)} & \textbf{\(N\)} & \textbf{\(I\)} & \(RMSE_u\) & \(MAE_u\) & \(MAPE_u\) & \(J^*_{\%error}\) & \(RMSE_u\) & \(MAE_u\) & \(MAPE_u\) & \(J^*_{\%error}\) \\
        \hline
        \hline
        \multirow{6}{*}{\rotatebox{90}{\textbf{1}\quad \textit{NN (Indirect)}}} 
        &\textbf{1} & 1 & - & - & 2 & 1.50E-02 & 1.03E-02 & 17.81 & 0.55 & 1.33E-02 & 9.40E-03 & 25.02 & 0.61 \\ \cline{7-14}
        &\textbf{2} & 1 & - & - & 6 & 7.09E-05 & 5.27E-05 & 0.10 & 0.02 & 5.91E-05 & 4.68E-05 & 0.17 & 0.03 \\ \cline{3-14}
        &\textbf{3} & 2 & - & - & 2 & 8.90E-03 & 5.50E-03 & 86.73 & 1.08E+03 & 8.00E-03 & 5.10E-03 & 122.79 & 1.44E+03 \\ \cline{7-14}
        &\textbf{4} & 2 & - & - & 6 & 1.95E-04 & 1.31E-04 & 1.75 & 2.35 & 1.69E-04 & 1.20E-04 & 2.24 & 1.97 \\ \cline{3-14}
        &\textbf{5} & 3 & - & - & 10 & 4.32E-01 & 2.43E-01 & 0.77 & 1.80 & 3.74E-01 & 2.25E-01 & 0.70 & 1.78 \\ \cline{7-14}
        &\textbf{6} & 3 & - & - & 30 & 3.87E-02 & 2.84E-02 & 0.09 & 0.09 & 3.65E-02 & 2.72E-02 & 0.09 & 0.09 \\
        
        \hline
        \hline
        \multirow{19}{*}{\rotatebox{90}{\textbf{2} \quad \textit{NN + Fourier Layer (Indirect)}}}
        &\textbf{7} & 1 & 4 & 4 & 2 & 3.70E-04 & 2.82E-04 & 0.39 & 0.06 & 3.48E-04 & 2.69E-04 & 0.50 & 0.07 \\ \cline{7-14}
        &\textbf{8}& 1 & 4 & 4 & 6 & 3.52E-04 & 2.56E-04 & 0.32 & 0.05 & 3.39E-04 & 2.50E-04 & 0.40 & 0.05 \\ \cline{7-14}
        &\textbf{9} & 1 & 5 & 5 & 2 & 4.71E-04 & 3.52E-04 & 1.02 & 0.07 & 4.29E-04 & 3.29E-04 & 1.50 & 0.08 \\ \cline{7-14}
        &\textbf{10} & 1 & 5 & 5 & 6 & 3.14E-04 & 2.34E-04 & 0.65 & 0.05 & 2.99E-04 & 2.28E-04 & 1.00 & 0.07 \\ \cline{7-14}
        &\textbf{11} & 1 & 6 & 4 & 2 & 3.55E-04 & 2.55E-04 & 0.36 & 0.04 & 3.46E-04 & 2.52E-04 & 0.48 & 0.04 \\ \cline{7-14}
        &\textbf{12} & 1 & 6 & 4 & 6 & 3.88E-04 & 2.84E-04 & 0.44 & 0.06 & 3.76E-04 & 2.79E-04 & 0.66 & 0.08 \\ \cline{3-14}
        &\textbf{13} & 2 & 4 & 4 & 2 & 2.20E-03 & 1.40E-03 & 38.28 & 6.33 & 1.80E-03 & 1.30E-03 & 51.98 & 10.51 \\ \cline{7-14}
        &\textbf{14} & 2 & 4 & 4 & 6 & 2.20E-03 & 1.40E-03 & 39.09 & 11.26 & 1.80E-03 & 1.30E-03 & 53.74 & 17.13 \\ \cline{7-14}
        &\textbf{15} & 2 & 5 & 5 & 2 & 6.57E-04 & 4.93E-04 & 10.26 & 3.90 & 5.58E-04 & 4.42E-04 & 11.52 & 5.25 \\ \cline{7-14}
        &\textbf{16} & 2 & 5 & 5 & 6 & 6.18E-04 & 4.51E-04 & 7.37 & 5.84 & 5.25E-04 & 3.97E-04 & 8.16 & 7.87 \\ \cline{7-14}
        &\textbf{17} & 2 & 6 & 4 & 2 & 2.20E-03 & 1.10E-03 & 37.70 & 7.48 & 1.70E-03 & 9.67E-04 & 50.59 & 12.35 \\ \cline{7-14}
        &\textbf{18} & 2 & 6 & 4 & 6 & 2.20E-03 & 1.10E-03 & 37.90 & 7.81 & 1.70E-03 & 9.79E-04 & 50.19 & 12.90 \\ \cline{7-14}
        &\textbf{19} & 2 & 8 & 8 & 2 & 2.63E-04 & 1.95E-04 &	3.83 &	2.32 & 2.28E-04 & 1.79E-04 &	4.94 & 3.91 \\ \cline{3-14}
        &\textbf{20} & 3 & 4 & 4 & 3 & 1.97E-01 & 1.34E-01 & 0.47 & 0.12 & 1.58E-01 & 1.20E-01 & 0.43 & 0.07 \\ \cline{7-14}
        &\textbf{21} & 3 & 4 & 4 & 6 & 1.72E-01 & 1.17E-01 & 0.41 & 0.10 & 1.40E-01 & 1.05E-01 & 0.37 & 0.06 \\ \cline{7-14}
        &\textbf{22} & 3 & 5 & 5 & 3 & 2.87E-01 & 1.93E-01 & 0.65 & 0.14 & 2.34E-01 & 1.76E-01 & 0.59 & 0.13 \\ \cline{7-14}
        &\textbf{23} & 3 & 5 & 5 & 6 & 5.87E-02 & 4.10E-02 & 0.14 & 0.06 & 4.92E-02 & 3.67E-02 & 0.13 & 0.06 \\ \cline{7-14}
        &\textbf{24} & 3 & 6 & 4 & 3 & 2.19E-01 & 1.23E-01 & 0.39 & 0.21 & 1.61E-01 & 1.06E-01 & 0.34 & 0.24 \\ \cline{7-14}
        &\textbf{25} & 3 & 6 & 4 & 6 & 2.17E-01 & 1.28E-01 & 0.44 & 0.17 & 1.63E-01 & 1.11E-01 & 0.39 & 0.12 \\
        \hline
        \hline
        \multirow{3}{*}{\rotatebox{90}{\parbox{1.5cm}{\textbf{3}\; \textit{Direct}}}}
        &\textbf{26} & 3 & 4 & 4 & - & 3.50E+00 & 2.17E+00 & 6.13 & 7.41 & 3.23E+00 & 2.08E+00 & 5.95 & 8.67 \\ \cline{7-14}
        &\textbf{27} & 3 & 5 & 5 & - & 2.92E+00 & 1.87E+00 & 5.48 & 3.59 & 2.65E+00 & 1.79E+00 & 5.32 & 4.42 \\ \cline{7-14}
        &\textbf{28} & 3 & 6 & 4 & - & 4.54E+00 & 2.85E+00 & 8.27 & 13.82 & 4.22E+00 & 2.73E+00 & 7.96 & 15.36 \\
        \hline
    \end{tabular}}
\end{table}

\clearpage

\bibliographystyle{apsrev4-1}
\bibliography{References}

\end{document}